\documentclass[titlepage,draft,12pt]{article} 
\usepackage{amsfonts,latexsym}
\usepackage[a4paper]{geometry}

\date{}

\newcommand{\ep}{\varepsilon}
\newcommand{\qed}{{\penalty 10000\mbox{$\quad\Box$}}}
\newcommand{\re}{\mathbb{R}}
\newcommand{\rebar}{\overline{\mathbb{R}}}
\newcommand{\n}{\mathbb{N}}

\newcommand{\uep}{u_{\ep}}
\newcommand{\uzep}{u_{0\ep}}
\newcommand{\vep}{v_{\ep}}
\newcommand{\Eep}{E_{\ep}}
\newcommand{\Fep}{E_{\ep}^{**}}
\newcommand{\phiep}{\varphi_{\ep}}
\newcommand{\phieps}{\varphi_{\ep}^{**}}
\newcommand{\Gammalim}{\Gamma\,\mbox{--}\!\lim}
\newcommand{\Gammaliminf}{\Gamma\,\mbox{--}\liminf}

\newtheorem{thm}{Theorem}[section]
\newtheorem{thmbibl}{Theorem}
\newtheorem{rmk}[thm]{Remark}
\newtheorem{prop}[thm]{Proposition}
\newtheorem{defn}[thm]{Definition}

\title{Passing to the limit in maximal slope curves: from a 
regularized Perona-Malik equation to the total variation flow}

\author{Maria Colombo\vspace{1ex}\\ 
{\normalsize Scuola Normale Superiore}\\
{\normalsize PISA (Italy)}\\
{\normalsize e-mail: \texttt{maria.colombo@sns.it}}
\and
Massimo Gobbino\vspace{1ex}\\ 
{\normalsize Universit\`a degli Studi di Pisa} \\
{\normalsize Dipartimento di Matematica Applicata ``Ulisse Dini''}\\ 
{\normalsize PISA (Italy)}\\  
{\normalsize e-mail: \texttt{m.gobbino@dma.unipi.it}}}

\begin{document}
\maketitle
\begin{abstract}
	We prove that solutions of a mildly regularized Perona-Malik 
	equation converge, in a slow time scale, to solutions of the 
	total variation flow. The convergence result is global-in-time, 
	and holds true in any space dimension.
	
	The proof is based on the general principle that ``the limit of
	gradient-flows is the gradient-flow of the limit''.  To this
	end, we exploit a general result relating the Gamma-limit of a
	sequence of functionals to the limit of the corresponding maximal
	slope curves.
	
\vspace{1cm}

\noindent{\bf Mathematics Subject Classification 2000 (MSC2000):}
35K55, 35B40, 35K90.

\vspace{1cm} 

\end{abstract}

 
\section{Introduction}

The Perona-Malik equation 
\begin{equation}
	u_{t}=\mathrm{div}\left(\frac{\nabla u}{1+|\nabla
	u|^{2}}\right)
	\quad\quad
	(x,t)\in\Omega\times(0,+\infty)
	\label{eqn:PM}
\end{equation}
(where $\Omega\subseteq\re^{n}$ is a bounded open set) is arguably the
most celebrated example of forward-backward diffusion process.  It was
introduced by P.~Perona and J.~Malik~\cite{PM} in the context of image
denoising.  It is formally the gradient-flow of the functional
\begin{equation}
	PM(u):=\frac{1}{2}\int_{\Omega}\log\left(1+|\nabla
	u(x)|^{2}\right)\,dx.
	\label{defn:PM-funct}
\end{equation}

The forward-backward nature of (\ref{eqn:PM}) depends on the
convex-concave behavior of the integrand in~(\ref{defn:PM-funct}).
Equation (\ref{eqn:PM}) has generated a considerable literature (see
\cite{amann,BFG,BF,BNP,BNP3,BNPT1,Z2,CLMC,nes,E1,E2,E3,GG,GG-tams,GG-cpdes,GG-ifb,G-PMEntire,
guidotti,gl,Z3,Z1}), focussed both on numerical and on analytical
aspects.  The big challenge is to reconcile the empirical and
practical efficacy of the method, supported by several numerical
observations, with the expected analytical ill-posedness of a backward
parabolic equation.  A satisfactory theory would represent a solution
to the Perona-Malik paradox, as named after~\cite{kich}, but it still
seems to be out of reach.

In a recent paper, P.\ Guidotti~\cite{patrick-DM} introduced a mild 
regularization of (\ref{eqn:PM}). He considered the family of 
functionals
$$ PM_{\delta}(u):=\frac{1}{2}\int_{\Omega}\log\left(1+|\nabla
	u(x)|^{2}\right)\,dx+ \delta
	\int_{\Omega}|\nabla u(x)|^{2}\,dx,$$
where $\delta>0$ is a parameter.  The integrand remains nonconvex, at
least in the interesting cases where $\delta$ is small, but now it is
convex-concave-convex, and it grows quadratically at infinity.  This
is enough to guarantee that the corresponding gradient-flow equation,
which simply reads as
\begin{equation}
	u_{t}=\mathrm{div}\left(\frac{\nabla u}{1+|\nabla
	u|^{2}}\right)+\delta\,\Delta u 
	\quad\quad
	(x,t)\in\Omega\times(0,+\infty),
	\label{eqn:PM-delta}
\end{equation}
has a unique global-in-time solution in the sense of Young measures,
according to the theory developed in \cite{KP,demoulini}.

Several qualitative properties of these solutions are reported in 
\cite{patrick-DM}. In particular, it seems that the well-known 
staircasing effect of the original Perona-Malik equation is now 
replaced by a ``ramping'' effect, namely the tendency of solutions to 
alternate flat plateaus and bounded growth regions in a piecewise 
fashion (see all figures in~\cite{patrick-DM}).

In \cite{patrick-DM} it is also observed that the limit of
approximated solutions as $\delta\to 0^{+}$ is the trivial stationary
solution frozen in the initial condition, and ``thus the only way to
produce a meaningful limit would involve a rescaling of time in the
process''.

In this paper we follow this path.  We take a family of solutions
$v_{\delta}(t)$ of the regularized model (with initial conditions
$v_{\delta}(0)$ converging to some $u_{0}$), and then we speed up the
evolution by considering the family of functions
\begin{equation}
	\uep(t):=v_{4^{-1}\ep^2|\log\ep|^{}}\left(\frac{t}{\ep|\log\ep|}\right)
	\quad\quad
	t\geq 0.
	\label{defn:uep}
\end{equation}

This happens to be the right rescaling factor, in the sense that 
$\uep(t)$ uniformly converges to a nontrivial limit $u(t)$ as $\ep\to 
0^{+}$. Moreover $u(t)$ turns out to be the solution of the total 
variation flow (see~\cite{TV-nonno,TV-esistenza})
\begin{equation}
	u_{t}=\mathrm{div}\left(\frac{\nabla u}{|\nabla u|}\right)
	\quad\quad
	(x,t)\in\Omega\times(0,+\infty),
	\label{eqn:TVF}
\end{equation}
with initial condition $u_{0}$.  Any slower time rescaling leads in
the limit to the stationary solution $u(t)\equiv u_{0}$, and any
faster time rescaling produces a limit which for every $t>0$
coincides with the constant function equal to the average of $u_{0}$
in $\Omega$.

All qualitative properties of the total variation flow
(see~\cite{TV-property}) are consistent with the numerical experiments
presented in~\cite{patrick-DM}.  Therefore, our convergence result
provides a rigorous justification of all these properties.

The proof of the convergence result involves three main steps.  First
of all, we interpret $\uep(t)$ as a gradient-flow.  As it comes
from~\cite{patrick-DM}, $\uep(t)$ is the solution in the sense of
Young measures of the forward-backward equation associated to the
formal gradient-flow of the nonconvex energy
\begin{equation}
	\Eep(u):=\int_{\Omega}
	\phiep\left(\strut|\nabla u(x)|\right)\,dx,
	\label{defn:Eep}
\end{equation}
where
$$\phiep(\sigma):=\frac{1}{2\ep|\log\ep|}
\log\left(1+\sigma^{2}\right)+\frac \ep 4\,\sigma^{2}.$$

The key point is that $\uep(t)$ is also the gradient-flow of 
the relaxed (convex) energy
\begin{equation}
	\Fep(u):=\int_{\Omega}\phieps\left(\strut|\nabla u(x)|\right)\,dx,
	\label{defn:Fep}
\end{equation}
where $\phieps$ is the convexification of $\phiep$.

In other words, as a result of the first step we can forget about
Young measures and forward-backward equations, and think of $\uep(t)$
as the solution of a degenerate forward parabolic equation, or better
as the gradient-flow of a convex (although not strictly convex)
functional.  We stress that this is a general fact.  Solutions
provided by the theory developed in~\cite{KP,demoulini} always
coincide with gradient-flows of the corresponding convexified (or
relaxed) energies.  Quite surprisingly, this has never been observed
in the literature up to our knowledge.

In the second step we compute the Gamma-limit of the energies
(Theorem~\ref{thm:Gamma-conv}), and we discover that
\begin{equation}
	\Gammalim_{\ep\to 0^{+}}\Eep(u)=
	\Gammalim_{\ep\to 0^{+}}\Fep(u)=
	TV(u),
	\label{th:Gamma-lim}
\end{equation}
where $TV(u)$ denotes the total variation of $u$.  Since the total
variation flow is the gradient-flow of $TV(u)$, our convergence result
is now equivalent to say that the limit of gradient-flows is the
gradient-flow of the Gamma-limit.

This is the content of the third step, where we deduce it from a
general result (Theorem~\ref{thm:main-cmp}) which we state and prove
in the abstract setting of \emph{maximal slope curves} in metric
spaces (see Section~\ref{sec:CMP} for further details and references).
We believe that the scope of the general result goes far beyond this
simple application.

We conclude by discussing how this problem relates to recent
investigations about the ``slow time'' behavior of approximations of
the Perona-Malik equation.  The starting point is the observation that
the evolution, despite of the different approximation methods, seems to
exhibit always three different time scales, named ``fast time'',
``standard time'', and ``slow time'' according to~\cite{BF}.

In a fast time of order $o(1)$ solutions develop microstructures in
the concave regime (staircasing).  In a time scale of order $O(1)$
(standard time) solutions behave as expected in the original model,
with a smoothing effect in the concave regime, and sharpening of
regions where the gradient is large.  At the ``end'' of standard time,
solutions have an almost piecewise constant structure, and this is consistent
with the intuitive idea that piecewise constant functions are
stationary points of $PM(u)$.  On the other hand, only constant
functions (and not piecewise constant functions) are stationary points
of the usual approximating models.  As a consequence, approximating
solutions exhibit a transition from a piecewise constant structure to
a constant value (equal to the average of $u_{0}$).  The transition
turns out to be very slow because there is almost no energy left.

This gives rise to the so called ``slow time'' motion, in which the
plateaus of the piecewise constant function move in the vertical
direction, with jump points which remain fixed in space. Although the 
existence of this phase seems to be independent of the approximation 
method, the law of the vertical motion does depend on it. 

All previous results on this problem are limited to the one
dimensional case.  Let $u$ be a piecewise constant function defined in
an interval, and let $S_{u}$ be the (finite or countable) set of its
jump points.  Let $J_{x}$ denote the jump height in a point $x\in
S_{u}$, and let us consider the following energies
$$H_{\alpha}(u):=\sum_{x\in S_{u}}|J_{x}|^{\alpha} \quad\mbox{(with
$\alpha\in(0,1]$)}, 
\quad\quad
H_{0}(u):=\sum_{x\in S_{u}}\log|J_{x}|.$$

In the case of a fourth order regularization of (\ref{eqn:PM}),
corresponding to adding a vanishing second order term to
(\ref{defn:PM-funct}), G.~Bellettini and A.~Fusco~\cite{BF}
conjectured that the vertical motion is governed by the gradient-flow
of $H_{1/2}(u)$.  They supported their conjecture by proving the
corresponding Gamma-limit result for the energies.  The missing step
is a rigorous proof that also in that case the limit of gradient-flows
is the gradient-flow of the limit.  Unfortunately
Theorem~\ref{thm:main-cmp} does not apply to their functionals.

In the case of the semidiscrete scheme (see~\cite{BNP3,nes}) the
vertical motion is governed by the gradient-flow of $H_{0}(u)$, which
in a certain sense represents a limit case.  The proof given
in~\cite{nes} exploits a variant of Theorem~\ref{thm:main-cmp},
complicated by the fact that $H_{0}(u)$ is not bounded from below.

What we show in this paper is that in the model proposed
in~\cite{patrick-DM} the slow-time vertical motion is governed by the
gradient-flow of $H_{1}(u)$.  This is the opposite limit case, and
phenomena are completely different.  The good news is that the limit
energy is convex.  This simplifies the analysis, which here can be
carried out in any space dimension, and delivers a well known limit
problem.  On the other hand, the relaxation of $H_{1}(u)$ is $TV(u)$,
hence it is finite in the whole space of bounded variation functions.
As a consequence, in this case the slow time motion is not limited to
piecewise constant functions.  This is hardly surprising after
reminding that in this model the motion in standard time is trivial.

This paper is organized as follows.  In Section~\ref{sec:statements}
we fix notations and we state our convergence results.
Section~\ref{sec:CMP} is devoted to limits of maximal slope curves in
metric spaces.  In Section~\ref{sec:proofs} we prove our main result.

\setcounter{equation}{0}
\section{Notations and statements}\label{sec:statements}

Let $n$ be a positive integer, and let $\Omega\subseteq\re^{n}$ be an
open set, which for simplicity we assume to be bounded and an
extension domain (see Definition~3.20 in~\cite{AFP}, satisfied by all
bounded open sets with Lipschitz boundary).  The more general ambient
space we consider is $L^{2}(\Omega)$.  We write
$\|u\|_{L^{p}(\Omega)}$, or simply $\|u\|_{p}$, to denote the $p$-norm
(with $p\in[1,+\infty]$) of a function $u\in L^{p}(\Omega)$.  All the
energies we consider, and in particular $\Eep(u)$ and $\Fep(u)$, are
always thought as defined in the whole space $L^{2}(\Omega)$ by
setting them equal to $+\infty$ outside their natural domain.  We
write $BV(\Omega)$ to denote the space of all functions $u\in
L^{2}(\Omega)$ with finite total variation $TV(u)$.  Once again, we
think of $TV(u)$ as defined for every $u\in L^{2}(\Omega)$, with
$TV(u)<+\infty$ if and only if $u\in BV(\Omega)$.

For every $\delta>0$ we consider equation (\ref{eqn:PM-delta}), with
Neumann boundary conditions, and an initial datum.  It has been shown
in~\cite{patrick-DM} that solutions $v_{\delta}(t)$ exist in a suitable weak
sense.  For every $\ep\in(0,1)$, we define $\uep(t)$ by rescaling
$v_{\delta}(t)$ according to (\ref{defn:uep}).  It turns out that
$\uep(t)$ is a solution in the same weak sense of equation
\begin{equation}
	{\uep}_{t}=\frac{1}{\ep|\log\ep|}\,\mathrm{div}\left(
	\frac{\nabla \uep}{1+|\nabla
	\uep|^{2}}\right)+\frac \ep 4\, \Delta \uep 
	\quad\quad
	(x,t)\in\Omega\times(0,+\infty),
	\label{eqn:PM-ep}
\end{equation}
with Neumann boundary conditions
\begin{equation}
	\frac{\partial\uep}{\partial n}(x,t)=0
	\quad\quad
	(x,t)\in\partial\Omega\times(0,+\infty),
	\label{eqn:PM-NBC}
\end{equation}
and initial condition
\begin{equation}
	\uep(x,0)=\uzep(x)
	\quad\quad
	x\in\Omega.
	\label{eqn:uep-data}
\end{equation}

In the following result we collect properties of $\uep(t)$.

\begin{thmbibl}[Properties of rescaled approximating solutions]\label{thm:uep}
	
	Let $n$ be a positive integer, let $\Omega\subseteq\re^{n}$ be a
	bounded extension domain, let $\ep\in (0,1)$, and let $\uzep\in
	L^{2}(\Omega)$. Let $\Fep$ be the functional defined in 
	(\ref{defn:Fep}). 
	
	Then the following properties hold true.
	\begin{enumerate}
		\renewcommand{\labelenumi}{(\arabic{enumi})} 
		\item
		\label{stat:yms} \emph{(Weak Young measure solution and
		regularity)} There exists a unique function $\uep(t)$ and a
		(not necessarily unique) gradient Young measure $\nu_{\ep}$ in
		$\Omega\times[0,+\infty)$ such that the pair
		$(\uep,\nu_{\ep})$ is a weak Young measure valued solution of
		problem (\ref{eqn:PM-ep}) through (\ref{eqn:uep-data}) in the
		sense of~\cite{KP,demoulini}.

		Moreover, we have that
		$$\uep\in C^{0}\left([0,+\infty);L^{2}(\Omega)\right)\cap
		C^{1}\left((0,+\infty);L^{2}(\Omega)\right),$$
		and for every $t>0$ we have that $\uep(t)$ is regular enough
		so that the right-hand side of (\ref{eqn:fd-pde}) lies in
		$L^{2}(\Omega)$.
	
		\item \emph{(Degenerate forward parabolic equation)} The
		function $\uep$ of statement~(\ref{stat:yms}) is the unique
		solution in $\Omega\times(0,+\infty)$ of the partial
		differential equation
		\begin{equation}
			{\uep}_{t}=-\nabla\Fep(\uep)=
			\mathrm{div}\left[(\phieps)'
			\left(\strut|\nabla\uep|\right)
			\frac{\nabla\uep}{|\nabla\uep|}\right],
			\label{eqn:fd-pde}
		\end{equation}
		with Neumann boundary conditions (\ref{eqn:PM-NBC}), and
		initial condition (\ref{eqn:uep-data}).

		\item \emph{(Gradient-flow integral inequality)} The 
		function $\uep$ of statement~(\ref{stat:yms}) 
		is the unique function satisfying (\ref{eqn:uep-data}) and 
		the inequality
		\begin{equation}
			\Fep(\uep(s))-\Fep(\uep(t))\geq
			\frac{1}{2}\int_{s}^{t}\|\uep'(\tau)\|_{2}^{2}\,d\tau+
			\frac{1}{2}\int_{s}^{t}\|\nabla\Fep(\uep(\tau))\|_{2}^{2}\,d\tau
			\label{hp:uep-gf}
		\end{equation}
		for every $0\leq s\leq t$.

		\item \label{stat:Linfty}\emph{($L^{p}$ estimate)}
		If $\uzep\in L^{p}(\Omega)$ for some $p\in[1,+\infty]$, then
		$\uep(t)\in L^{p}(\Omega)$ for every $t\geq 0$, and the function
		$t\to\|\uep(t)\|_{L^{p}(\Omega)}$ is nonincreasing.
	\end{enumerate}	

\end{thmbibl}

Theorem~\ref{thm:uep} above shows that the function $\uep(t)$ can be
characterized in at least three different ways, either as the solution
in the sense of Young measures of a forward-backward equation, or as
the solution of a forward degenerate parabolic equation (for example
in the sense of~\cite{brezis}), or as a maximal slope curve
(gradient-flow inequalities).  What we need in this paper is only the
last one.  Solutions generate a contraction semigroup in
$L^{2}(\Omega)$.

As far as we know this equivalence has never been stated explicitly. 
Nevertheless, it follows from some general facts which nowadays are 
quite well known, and which we now recall briefly.

First of all, the three approaches lead to a \emph{unique} solution. 
Uniqueness follows in the first case from the strong requirements 
imposed on the structure of the corresponding Young measure, in the 
second case from the contraction property of the semigroup generated 
by a forward parabolic equation, in the third case from the convexity 
of the energy $\Fep(u)$.

Secondly, in all three approaches the solution is usually obtained (or
at least it can be obtained) as the limit of approximated solutions
constructed via an iterated minimization process, known as minimizing
movement (see \cite{dg-min}).  As already observed
in~\cite{KP,demoulini,patrick-DM}, the minimization procedure gives
the same result when applied to $\Eep(u)$ or $\Fep(u)$.

In conclusion, all three approaches define a unique solution through
an analogous procedure, hence the solution is the same.

In this paper we are interested in the behavior of $\uep(t)$ as 
$\ep\to 0^{+}$. Following the gradient-flow approach, the first thing 
to do is understanding the limit behavior of the energies. This is 
the content of next result.

\begin{thm}[Gamma-convergence and compactness]\label{thm:Gamma-conv}
	Let $n$ be a positive integer, and let $\Omega\subseteq\re^{n}$ be
	a bounded extension domain.  For every $\ep>0$, let
	$\Eep:L^{2}(\Omega)\to[0,+\infty]$ and
	$\Fep:L^{2}(\Omega)\to[0,+\infty]$ be defined as in
	(\ref{defn:Eep}) and (\ref{defn:Fep}), respectively, if $u\in
	H^{1}(\Omega)$, and $+\infty$ otherwise.
	
	Then we have the following conclusions.
	\begin{enumerate}
		\renewcommand{\labelenumi}{(\arabic{enumi})}
			
		\item \emph{(Gamma-convergence)} We have that
		(\ref{th:Gamma-lim}) holds true with respect to the topology
		of $L^{2}(\Omega)$.
		
		\item \emph{(Compactness)} Let
		$\{\uep\}_{\ep\in(0,1)}\subseteq L^{2}(\Omega)$ be a family 
		of functions such that
		\begin{equation}
			\sup_{\ep\in(0,1)}\left\{\strut\|\uep\|_{\infty}+
			\Fep(\uep)\right\}<+\infty.
			\label{hp:cpt}
		\end{equation}
		
		Then the family $\{\uep\}$ is relatively compact in
		$L^{2}(\Omega)$.
	\end{enumerate}
\end{thm}

The gradient-flow of the limit functional $TV(u)$ is the so called
\emph{total variation flow}, and the corresponding partial
differential equation is (\ref{eqn:TVF}).  The right-hand side of
(\ref{eqn:TVF}) needs to be interpreted in a suitable weak sense when
the gradient vanishes, and this happens in large regions because
solutions tend to develop flat plateaus.  A general existence and
uniqueness result was proved by F.\ Andreu, C.\ Ballester, V.\
Caselles, and J.\ M.\ Maz\'{o}n~\cite{TV-esistenza} (see
also~\cite{TV-property}) using the theory of accretive operators in
Banach spaces.  In~\cite{TV-esistenza} the operator in the right-hand
side of~(\ref{eqn:TVF}) is interpreted as the limit of the
$p$-Laplacian as $p\to 1^{+}$.

In this paper we limit ourselves to initial data in $L^{2}(\Omega)$,
in which case existence of a unique solution is provided also by the
theory of maximal monotone operators~\cite{brezis}, as explained
in~\cite{TV-nonno}.  In this context the right-hand side
of~(\ref{eqn:TVF}) is the subdifferential of the convex functional
$TV(u)$.  As in the case of approximating problems, this formulation
is equivalent to the gradient-flow integral inequality
	$$TV(u(s))-TV(u(t))\geq
	\frac{1}{2}\int_{s}^{t}\|u'(\tau)\|_{2}^{2}\,d\tau+
	\frac{1}{2}\int_{s}^{t}\|\nabla TV(u(\tau))\|_{2}^{2}\,d\tau$$
for every $0\leq s\leq t$, where $\|\nabla TV(u(\tau))\|_{2}$ is the
minimal norm of an element in the subdifferential of the functional
$TV$ in the point $u(\tau)$, which in turn coincides with the slope
of $TV$ in $u(\tau)$ as defined in Section~\ref{sec:CMP} in an
abstract metric setting.  This is the characterization of the total
variation flow which we need in this paper.

Our main result is the convergence of $\uep(t)$ to the solution of the
total variation flow with the same boundary conditions, and the limit
initial datum.  We point out that we do not assume initial data to be
a recovery sequence, and we do not ask their energy to be bounded.

\begin{thm}[Global-in-time convergence]\label{thm:main}
	Let $n$ be a positive integer, let $\Omega\subseteq\re^{n}$ be a
	bounded extension domain, let $u_{0}\in L^{2}(\Omega)$, and let
	$\{\uzep\}_{\ep\in(0,1)}\subseteq L^{2}(\Omega)$ be a family of
	functions such that
	\begin{equation}
		\lim_{\ep\to 0^{+}}\uzep=u_{0}
		\quad\quad
		\mbox{in }L^{2}(\Omega).
		\label{hp:data-conv}
	\end{equation}
	
	For every $\ep\in(0,1)$, let $\uep$ be the solution of the
	rescaled approximating problem with initial condition $\uzep$, in
	the sense of Theorem~\ref{thm:uep}.  Let $u(t)$ be the solution of
	the total variation flow with Neumann boundary conditions and
	initial datum $u_{0}$.
	
	Then we have that $\uep(t)\to u(t)$ in
	$C^{0}\left([0,+\infty);L^{2}(\Omega)\right)$, namely
	\begin{equation}
		\lim_{\ep\to 0^{+}}\sup_{t\geq 0}
		\|\uep(t)-u(t)\|_{L^{2}(\Omega)}=0.
		\label{th:main}
	\end{equation}
	
\end{thm}

We conclude this section with a heuristic argument which justifies 
the rescaling leading to (\ref{defn:Eep}).

Let us consider the one dimensional case where $\Omega=(-1,1)$, let
$J>0$, and let $v\in BV((-1,1))$ be the piecewise constant function
equal to $-J/2$ in $(-1,0)$, and equal to $J/2$ in $(0,1)$.  Let
$\eta>0$, and let $h\in H^{1}(\re)$ be a function such that
$h(x)= J/2$ for every $x\geq \eta$, $h(x)= -J/2$ for every $x\leq 
-\eta$, and  $h'(x)> 0$ for almost every $x\in (-\eta, \eta)$. 
For every $\ep>0$, let $\vep(x)$ be the approximation of $v$ defined as $\vep(x):=
h(x/\ep)$ for every $x\in (-1,1)$.

For every $\ep \leq \eta^{-1}$, plugging $\vep$ into (\ref{defn:Eep}),
with a variable change we obtain that
$$\Eep(\vep)=\frac{1}{2|\log\ep|}\int_{-\eta}^{\eta}
\log\left(1+\frac{1}{\ep^{2}}\left[h'(x)\right]^{2}\right)dx+
\frac 1 4 \int_{-\eta}^{\eta}\left[h'(x)\right]^{2}dx.$$

Letting $\ep\to 0^{+}$ we find that
\begin{equation}
	\liminf_{\ep\to 0^{+}}\Eep(\vep)\geq
	2\eta+\frac 1 4 \int_{-\eta}^{\eta}\left[h'(x)\right]^{2}dx.
	\label{heu:liminf}
\end{equation}

In order to estimate the right-hand side, we first minimize with
respect to $h$, and we discover that the optimal choice is the
function $h(x)$ defined as $J(2\eta)^{-1}x$ for every
$x\in(-\eta,\eta)$.  Then we compute the integral, and finally we
apply the inequality between arithmetic and geometric mean to deduce
that
\begin{equation}
	2\eta+\frac 1 4 \int_{-\eta}^{\eta}\left[h'(x)\right]^{2}dx\geq
	2\eta+\frac{J^{2}}{8\eta}\geq J.
	\label{heu:AM-GM}
\end{equation}

Estimates (\ref{heu:liminf}) and (\ref{heu:AM-GM}) suggest that the
cost of a jump could be the jump height, which leads to conjecture
that the Gamma-limit is $TV(u)$ for a general $u$ in any dimension.

\setcounter{equation}{0}
\section{Passing to the limit in maximal slope curves}\label{sec:CMP}

The abstract theory of gradient-flows in metric spaces was introduced
in~\cite{dgmt}, and then developed by the same authors and
collaborators in a series of papers (see~\cite{DMT,MST-SNS} and the
references quoted therein).  For a modern presentation we refer
to~\cite{AGS}.  Here we just recall some basic definitions.

Let $(X,d)$ be a metric space, let
$\rebar:=\re\cup\{-\infty,+\infty\}$ be the extended real line, and
let $F:X\to\rebar$ be any function.  The (descending) \emph{slope}
$|\nabla F|(x)$ of $F$ in $x$ is defined to be $+\infty$ if
$F(x)\not\in\re$, and otherwise
$$|\nabla F|(x):=\limsup_{y\to x}\frac{\max\{F(x)-F(y),0\}}{d(x,y)}
\in[0,+\infty].$$

For every $T>0$, the space $AC^{2}\left([0,T];X\right)$ is the set of 
all functions $v:[0,T]\to X$ for which there exists $g\in 
L^{2}((0,T))$ such that
\begin{equation}
	d(v(t),v(s))\leq\int_{s}^{t}g(\tau)\,d\tau
	\quad\quad
	\forall\, 0\leq s\leq t\leq T.
	\label{defn:metric-deriv}
\end{equation}

It can be seen that there exists a smallest function $g(t)$
satisfying (\ref{defn:metric-deriv}).  This function is called the
\emph{metric derivative} of $v$, and it is denoted by $|v'|(t)$.

A \emph{maximal slope curve} for $F$ in $[0,T]$ is a triple
$(u,\psi,E)$ where
\begin{itemize}
	\item  $u\in AC^{2}\left([0,T];X\right)$,

	\item  $\psi:[0,T]\to\re$ is a nonincreasing function  
	such that for every $0\leq s\leq t\leq T$ we have that
	\begin{equation}
		\psi(s)-\psi(t)\geq \frac 1 2
		\int_s^t |u'|^2(\tau)\,d\tau+ 
		\frac 1 2 \int_s^t  \left|\nabla F \right|^2
		(u(\tau))\,d\tau,
		\label{eqn:cmp-lim-psi}
	\end{equation}

	\item $E\subseteq [0,T]$ is a set with Lebesgue measure
	equal to 0 such that
	\begin{equation}
		\psi(t)=F(u(t))
		\quad\quad
		\forall t\in[0,T]\setminus E.
		\label{eqn:psi-ae}
	\end{equation}
\end{itemize}

To be more precise, the second integral in the right-hand side of
(\ref{eqn:cmp-lim-psi}) should be an upper integral, since at this level of
generality there is no reason for the function $t\to|\nabla F|(u(t))$
to be measurable.  On the other hand, it can be easily proved that
actually it is always true that $|u'|(t)=|\nabla F|(u(t))$ for almost
every $t\in[0,T]$, which implies the required measurability.

When $F$ is a $C^{1}$ function in a Hilbert space $X$, this weak 
formulation is equivalent to the classical one, namely to asking that 
$u'(t)=-\nabla F(u(t))$ for every $t\in[0,T]$. 

Besides generality, the advantage of this weak formulation is that
inequalities and integrals are more stable than equalities and
derivatives.  It follows that maximal slope curves exist under general
assumptions on $F$ (see Theorem 2.3.1 in~\cite{AGS}), and are quite
stable when passing to the limit, both with respect to initial
conditions, and with respect to functionals.  Results in this
direction are contained in~\cite{DMT} and~\cite{ss-cpam} in a Hilbert
setting, and in~\cite{ss-09} in a metric setting, but assuming that
initial data are a recovery sequence.

Here we state a quite general result, used in a special case also
in~\cite{nes}.

\begin{thm}[Limits of maximal slope curves]\label{thm:main-cmp}
	Let $X$ be a metric space, let $F: X \to\rebar$ be a
	function, and let $F_n: X \to\rebar$ be a sequence of
	functions. 
	
	Let us assume that for every $x\in X$, and every pair of sequences
	$\{n_{k}\}\subseteq\n$ and $\{x_{k}\}\subseteq X$, we have the
	implication
	\begin{equation}
		\fbox{$\begin{array}{r}
			n_{k}\to +\infty,\ x_{k}\to x  \\
			\noalign{\vspace{1ex}}
			\displaystyle{\sup_{k \in \n} 
			\left\{|F_{n_{k}}(x_{k})|+|\nabla F_{n_{k}}
			|(x_{k})\strut\right\} < +\infty}
		\end{array}$}
		\Longrightarrow
		\fbox{$\begin{array}{l}
			\displaystyle{\lim_{k\to +\infty}F_{n_{k}}(x_{k})=F(x).} \\
			\noalign{\vspace{1ex}}
			\displaystyle{\liminf_{k\to +\infty}|\nabla
			F_{n_{k}}|(x_{k}) \geq |\nabla F|(x).}
		\end{array}$}
		\label{hp:slope-energy}
	\end{equation}
	
	Let $T>0$, and for every $n\in \n$ let $(u_n,\psi_{n},E_{n})$ be
	a maximal slope curve for $F_n$ in $[0,T]$.  Let us assume that
	\begin{equation}
		\sup_{n\in\n}\;\max_{t\in[0,T]}|\psi_{n}(t)|
		< +\infty,
		\label{hp:sup-fn}
	\end{equation}
	and that $\{u_{n}\}$ has a pointwise limit, namely there
	exists $u:[0,T]\to X$ such that
	$$\lim_{n\to +\infty}u_{n}(t)=u(t)
	\quad\quad
	\forall t\in [0,T].$$
	
	Then there exist $\psi$ and $E$ such that $(u,\psi,E)$ is a
	maximal slope curve for $F$.
\end{thm}

\paragraph{\textmd{\emph{Proof}}}

From the definition of maximal slope curve, for every
$n\in\n$ we have that
\begin{equation}
	\psi_n(s)-\psi_n(t)\geq \frac 1 2 \int_s^t
	|u_n'|^2(\tau)\,d\tau+ \frac 1 2 \int_s^t \left|\nabla F_n \right|^2
	(u_n(\tau))\,d\tau
	\label{eqn:cmpn}
\end{equation}
for every $0\leq s\leq t\leq T$, and 
\begin{equation}
	\psi_{n}(t)=F_{n}(u_{n}(t))
	\quad\quad
	\forall t\in[0,T]\setminus E_{n}.
	\label{eqn:psi-ae-n}
\end{equation}

The functions $\psi_n(t)$ are nonincreasing, and equi-bounded because
of (\ref{hp:sup-fn}).  Thus the usual compactness result for monotone
functions (known as Helly's Lemma, see for
example~\cite[Lemma~3.3.3]{AGS}) implies the existence of a
nonincreasing function $\psi:[0,T]\to [0,+\infty)$ such that (up to
subsequences, not relabeled)
\begin{equation}
	\lim_{n\to +\infty}\psi_n(t)=\psi(t)
	\quad\quad
	\forall t\in[0,T].
	\label{defn:psi}
\end{equation}

Let us consider now the right-hand side of
(\ref{eqn:cmpn}).
Setting $s=0$ and $t=T$, and using once more assumption 
(\ref{hp:sup-fn}), we obtain that
\begin{equation}
	\sup_{n\in\n}\int_{0}^{T}|u_n'|^2(\tau)\,d\tau<+\infty,
	\label{est:t1-deriv}
\end{equation}
\begin{equation}
	\sup_{n\in\n}\int_{0}^{T}\left|\nabla F_n \right|^2
	(u_n(\tau))\,d\tau<+\infty.
	\label{est:t1-slope}
\end{equation}

From (\ref{est:t1-deriv}) we easily deduce that $u\in
AC^{2}\left([0,T];X\right)$, and
\begin{equation}
	\liminf_{n\to +\infty}\int_s^t |u_n'|^2(\tau)\,d\tau\geq
	\int_s^t |u'|^2(\tau)\,d\tau
	\quad\quad
	\forall\,0\leq s\leq t\leq T.
	\label{est:liminf-deriv}
\end{equation}

From (\ref{est:t1-slope}) and Fatou's Lemma we obtain that
$$\int_{0}^{T_{}} \left(\liminf_{n\to +\infty}\left|\nabla F_n
\right|^2 (u_n(\tau))\right)\,d\tau\leq \liminf_{n\to
+\infty}\int_{0}^{T_{}} \left|\nabla F_n \right|^2
(u_n(\tau))\,d\tau<+\infty,$$
hence there exists a set $E' \subseteq[0,T]$, with Lebesgue 
measure equal to 0, such that
\begin{equation}
	\liminf_{n\to +\infty}\left|\nabla F_n \right|
	(u_n(t))<+\infty
	\quad\quad
	\forall t\in[0,T]\setminus E'.
	\label{est:slliminf}
\end{equation}

Let us introduce the set
$$E := E' \cup \left( \bigcup_{n\in\n} E_n\right),$$
which has clearly Lebesgue measure equal to 0.  Let us consider any
$t\in[0,T]\setminus E$. Since $t\not\in E_{n}$, from 
(\ref{eqn:psi-ae-n}) and (\ref{hp:sup-fn}) we have that
\begin{equation}
	\sup_{n\in\n}\left|F_{n}(u_{n}(t)) \right| =
	\sup_{n\in\n}|\psi_{n}(t)|<+\infty.
	\label{est:bounded-energy}
\end{equation}

Moreover, due to (\ref{est:slliminf}) there exists a ($t$-dependent)
sequence $n_{k}\to +\infty$ such that
\begin{equation}
	\lim_{k\to +\infty}|\nabla F_{n_{k}}|(u_{n_{k}}(t))=
	\liminf_{n\to +\infty}|\nabla F_{n}|(u_{n}(t))<+\infty.
	\label{est:bounded-slope}
\end{equation}

Thanks to (\ref{est:bounded-energy}) and (\ref{est:bounded-slope}),
sequences $\{n_{k}\}$ and $\{u_{n_{k}}(t)\}$ satisfy the assumptions
in the left-hand side of (\ref{hp:slope-energy}).  If follows that for
every $t\in[0,T]\setminus E$ we have that
$$\psi(t)=\lim_{n\to +\infty}\psi_{n}(t)=\lim_{k\to +\infty}\psi_{n_k}(t)=
\lim_{k\to +\infty}F_{n_{k}}(u_{n_{k}}(t))=
F(u(t)),$$
which proves (\ref{eqn:psi-ae}), and
$$|\nabla F|(u(t))\leq \liminf_{k\to +\infty}|\nabla
F_{n_{k}}|(u_{n_{k}}(t))=\! 
\lim_{k\to +\infty}|\nabla F_{n_{k}}|(u_{n_{k}}(t))=
\liminf_{n\to +\infty}|\nabla
F_{n}|(u_{n}(t)),$$
so that one more application of Fatou's Lemma gives that
\begin{equation}
	\int_s^t \left|\nabla F \right|^2 (u(\tau))\,d\tau\leq
	\liminf_{n\to +\infty}\int_s^t\left|\nabla F_n \right|^2
	(u_n(\tau))\,d\tau
	\quad\quad
	\forall\,0\leq s\leq t\leq T.
	\label{est:liminf-slope}
\end{equation}

We can now take the $\liminf$ of both sides of (\ref{eqn:cmpn}).
Thanks to (\ref{defn:psi}), (\ref{est:liminf-deriv}), and
(\ref{est:liminf-slope}) we obtain that (\ref{eqn:cmp-lim-psi}) holds 
true. This completes the proof that $(u,\psi,E)$ is a maximal slope 
curve for $F$.\qed

\medskip

Thanks to Theorem~\ref{thm:main-cmp} above, any convergence result for
maximal slope curves is reduced to verifying three assumptions.  The
first one is the existence of a pointwise limit, namely a compactness
result.  The second one is estimate (\ref{hp:sup-fn}), which in
general follows from suitable assumptions on the sequence of initial
data and some boundedness from below of the functionals.  The third and
more important assumption is (\ref{hp:slope-energy}).  In the last
part of this section we show that (\ref{hp:slope-energy}) follows from
Gamma-convergence in a class of functionals which contains all convex
functionals in Banach spaces.

\begin{defn}[Slope Cone Property]\label{defn:scp}
	\begin{em}
		Let $(X,d)$ be a metric space.  A function $F: X \to\rebar$
		satisfies the \emph{Slope Cone Property} if
		$$F(y) \geq F(x) - |\nabla F |(x)\cdot d(x,y)$$
		for every $y\in X$ and every $x\in X$ such that $F(x)\in\re$ 
		and $|\nabla F |(x)<+\infty$.
	\end{em}
\end{defn}

\begin{rmk}\label{rmk:convex->scp}
	\begin{em}
		Let $X$ be a Banach space.  Then every \emph{convex} function
		$F: X \to [0, +\infty]$ fulfils the Slope Cone Property.  If
		in addition $F$ is lower semicontinuous, then its slope
		$|\nabla F|(x)$ is lower semicontinuous.
		
		On the other hand, also in Banach spaces there do exist
		nonconvex functions satisfying the Slope Cone Property.  An
		example is $F(x)=-|x|$ in $\re$.
	\end{em}
\end{rmk}

\begin{prop}\label{prop:scp}
	Let $X$ be a metric space,  and let $F_n: X \to\rebar$ be a sequence of
	functions satisfying the Slope Cone Property. Let us assume that 
	there exists
	\begin{equation}
		F(x):=\Gammalim_{n\to +\infty} F_n(x).
		\label{hp:gamma-conv}
	\end{equation}
	
	Then the sequence $\{F_{n}\}$ satisfies assumption
	(\ref{hp:slope-energy}) of Theorem~\ref{thm:main-cmp}.
\end{prop}

\paragraph{\textmd{\emph{Proof}}}

Let $n_{k}\to +\infty$ and $x_{k}\to x$ be two sequences as in
assumption (\ref{hp:slope-energy}).  Let $M$ be the supremum in the
left-hand side of (\ref{hp:slope-energy}).  Let $z_{n}\to x$ be a
recovery sequence for $x$, namely a sequence such that
$F_{n}(z_{n})\to F(x)$.

From the Slope Cone Property, and the uniform bound on slopes, it
follows that 
$$F_{n_{k}}(z_{n_{k}}) \geq F_{n_{k}}(x_{k}) - |\nabla F_{n_{k}}
|(x_{k})\cdot d(x_{k},z_{n_{k}}) \geq F_{n_{k}}(x_{k}) - M \,
d(x_{k},z_{n_{k}}).$$

Taking the $\limsup$ of both sides we obtain that
$$F(x) = \limsup_{k\to+\infty} F_{n_{k}}(z_{n_{k}}) \geq \limsup_{k\to+\infty}
F_{n_{k}}(x_k).$$

The opposite inequality with the $\liminf$ follows from assumption
(\ref{hp:gamma-conv}).  This proves the first limit in the right-hand
side of (\ref{hp:slope-energy}).

Let us prove now the $\liminf$ inequality for slopes.  Let $L$ denote
the $\liminf$ in the left-hand side of (\ref{hp:slope-energy}).  Let
us take any $y\in X$, and let $y_{n}\to y$ be a corresponding recovery
sequence.  Due to the Slope Cone Property we have that
$$F_{n_{k}}(y_{n_{k}})\geq F_{n_{k}}(x_{k})-
|\nabla F_{n_{k}}|(x_{k})\cdot d(x_{k},y_{n_{k}}).$$

We already proved that the first term in the right-hand side tends to
$F(x)$.  Therefore, taking the $\limsup$ of both sides we obtain that
$F(y)\geq F(x)-L\,d(x,y)$.  Since $y$ is arbitrary, this easily
implies that $L\geq|\nabla F|(x)$, which completes the proof.\qed
\medskip

We conclude by mentioning a straightforward extension of 
Definition~\ref{defn:scp} and Proposition~\ref{prop:scp} in the same 
spirit of \cite{DMT, MST-SNS}.

\begin{rmk}
	\begin{em}
		One can weaken Definition~\ref{defn:scp} by asking that
		$F:X\to\re\cup\{+\infty\}$ 
		(so we exclude the value $-\infty$), 
		and there exists a continuous function
		$\Phi:X^{2}\times\re^{3}\to\re$ such that $\Phi(x,x,u,v,w)=0$
		for every $(x,u,v,w)\in X\times\re^{3}$, and
		$$F(y) \geq F(x) - |\nabla F |(x)\cdot d(x,y)-\Phi\left(\strut
		x,y,F(x),F(y),|\nabla F|(x)\right)\cdot d(x,y)$$
		for every $(x,y)\in X^{2}$ such that $F(x)\in\re$, 
		$F(y)\in\re$, and $|\nabla F|(x)\in\re$.
		
		We call this property ``$\Phi$ Slope Cone Property''.  In
		a Banach space it is fulfilled, for example, by the sum of a
		$C^{1}$ function and a convex function.
		
		It can be easily proved that Proposition~\ref{prop:scp} holds
		true also if all functions $F_{n}$ satisfy the $\Phi$ Slope
		Cone Property with respect to the same function $\Phi$.
	\end{em}
\end{rmk}

\setcounter{equation}{0}
\section{Proofs}\label{sec:proofs}

\subsection{Gamma-convergence and compactness}

The first equality in (\ref{th:Gamma-lim}) is a general property of
Gamma-convergence.  So we can concentrate on the second one.  

A standard approach, suggested by the heuristic argument at the end of
Section~\ref{sec:statements}, involves a blow up argument in order to
reduce the Gamma-liminf inequality to minimizing the right-hand side
of (\ref{heu:liminf}) with respect to $\eta$ and $h$, and a proof of
the Gamma-limsup inequality via the density of piecewise constant
functions with smooth level sets, for which a recovery sequence can be
constructed by adapting the optimal profile in the direction
orthogonal to level sets.

In both cases we follow a different and more elementary approach, 
which exploits that our functionals depend only on the gradient.

\paragraph{\emph{\textmd{Gamma-liminf inequality}}}

We claim that for every $a\in(0,1)$ and $b\in(0,1)$ there exists 
$\ep_{1}\in(0,1)$ such that 
\begin{equation}
	\phieps(\sigma)\geq a|\sigma|-b
	\quad\quad
	\forall\sigma\in\re,\ \forall\ep\in(0,\ep_{1}).
	\label{th:phiep>}
\end{equation}

If we prove this claim, then for every $u\in L^{2}(\Omega)$ and every
$\ep\in(0,\ep_{1})$ we have that 
$$\Fep(u)=\int_\Omega \phieps\left(\strut|\nabla u(x)|\right)\,dx \geq 
\int_\Omega
\left(\strut a|\nabla u(x)|-b\right)dx = a\, TV(u)-b|\Omega|.$$

Since the functional $TV(u)$ is lower 
semicontinuous, this proves that
	$$\Gammaliminf_{\ep \to 0^{+}}\Fep(u) \geq  a\,TV(u)-b|\Omega|.$$

Letting $a\to 1^{-}$ and $b\to 0^{+}$, we obtain the required
inequality.

So we are left to prove (\ref{th:phiep>}). Since the right-hand side 
is convex, it is enough to prove (\ref{th:phiep>}) with $\phiep$ 
instead of $\phieps$. Moreover, without loss of generality we can 
assume that $\sigma\geq 0$. Now we distinguish four cases.

If $\sigma\in[0,b]$, then we have that
$$\phiep(\sigma)\geq 0\geq ab-b\geq a \sigma-b.$$

If $\sigma\in[b,(e^{2}-1)^{1/2}]$, then we have that
$$\phiep(\sigma)-a\sigma\geq\frac{1}{2\ep|\log\ep|}\log(1+b^{2})-
a\sqrt{e^{2}-1}.$$

The right-hand side tends to $+\infty$ as $\ep\to 0^{+}$, so it is 
greater than $-b$ when $\ep$ is small enough.

If $\sigma\in[(e^{2}-1)^{1/2},(\ep|\log\ep|)^{-1}]$, then we have that
$$\phiep(\sigma)-a\sigma\geq\frac{1}{2\ep|\log\ep|}\log e^{2}-
\frac{a}{\ep|\log\ep|}=
\frac{1-a}{\ep|\log\ep|},$$
so that the conclusion follows as in the previous case.

Finally, when $\sigma\geq(\ep|\log\ep|)^{-1}$ we apply the inequality
between arithmetic mean and geometric mean, and we obtain that
$$\phiep(\sigma)\geq \frac{\log\sigma}{\ep|\log\ep|}+
\frac{\ep}{4}\,\sigma^{2}\geq\sigma\cdot
\left(\frac{\log\sigma}{|\log\ep|}\right)^{1/2}
\geq\sigma\cdot\left\{\frac{1}{|\log\ep|}\log\left(
\frac{1}{\ep|\log\ep|}\right)\right\}^{1/2}.$$

It is not difficult to see that the coefficient of $\sigma$ tends to 
1 as $\ep\to 0^{+}$, hence it is greater than $a$ when $\ep$ is small 
enough (in this point it is essential that $a<1$).

This completes the proof of (\ref{th:phiep>}), hence also of
the Gamma-liminf inequality.

\paragraph{\emph{\textmd{Gamma-limsup inequality}}}

By a classical density argument it is enough to find a recovery 
sequence for all functions $u\in C^{1}(\Omega)$ whose gradient is 
bounded in $\Omega$. To this end, it is enough to show that for any 
such function we have that
\begin{equation}
	\limsup_{\ep\to 0^{+}}\int_{\Omega}
	\phieps\left(\strut|\nabla u(x)|\right)\,dx
	\leq TV(u).
	\label{th:gamma-limsup}
\end{equation}

In turn, (\ref{th:gamma-limsup}) is proved if we show that for every 
$M>0$ we have that
	$$\phieps(\sigma)\leq a_{\ep}|\sigma|
	\quad\quad
	\forall\sigma\in[-M,M]$$
for a suitable coefficient $a_{\ep}$ which tends to $1$ as $\ep\to
0^{+}$.  Let us assume, without loss of generality, that
$\sigma\in[0,M]$.  Since $\phieps$ is the convexification of $\phiep$,
we can estimate $\phieps(\sigma)$ from above with the linear function
interpolating the values of $\phiep$ in $0$ and $2\ep^{-1}$.  As soon
as $M\geq 2\ep^{-1}$ we obtain that
$$\phieps(\sigma)\leq\frac \ep
2\phiep\left(\frac 2 \ep\right)\cdot\sigma \leq \frac 1 2
\left(\frac {\log(1+4\ep^{-2})}{2|\log \ep|}+ 1\right)\cdot\sigma 
\quad\quad
\forall\sigma\in [0,M].$$

As required, the coefficient of $\sigma$ in the right-hand side tends
to $1$.  This completes the proof of (\ref{th:gamma-limsup}).

\paragraph{\emph{\textmd{Compactness}}}

From (\ref{th:phiep>}) with $a=b=1/2$ we obtain that
$$\Fep(\uep) \geq \frac 1 2 \int_\Omega
\left(\strut|\nabla\uep(x)| - 1\right)dx = \frac 1 2 \,TV(\uep) - \frac 1 
2|\Omega|.$$

This estimate and assumption (\ref{hp:cpt}) yield a uniform bound on
the $L^{\infty}$-norm and on the total variation of $\uep$.  Thanks to
well known embedding theorems in $BV(\Omega)$ (see~\cite{AFP}, this is
the point where we need $\Omega$ to be an extension domain), this
implies that the family $\{\uep\}$ is relatively compact in
$L^{p}(\Omega)$ for every $p<+\infty$, and in particular in
$L^{2}(\Omega)$.\qed

\subsection{Convergence of approximating solutions}

In the first three paragraphs we prove the result with the further 
assumption that
\begin{equation}
	\sup_{\ep\in(0,1)}\left\{\strut\|\uzep\|_{\infty}+
	\Fep(\uzep)\right\}<+\infty.
	\label{hp:uzep}
\end{equation}

Then in the last paragraph we prove it for general data.

\paragraph{\emph{\textmd{Compactness on bounded time intervals}}}

We show that, for every $T>0$, the family $\{\uep\}$ is relatively
compact in $C^{0}\left([0,T];L^{2}(\Omega)\right)$.  This
follows from Ascoli's theorem provided that we show that solutions are
$1/2$-H\"{o}lder continuous with equi-bounded H\"{o}lder constants,
and that for every fixed $t\geq 0$ the family $\{\uep(t)\}\subseteq
L^{2}(\Omega)$ is relatively compact.

From (\ref{hp:uep-gf}) and H\"{o}lder's inequality we have 
that
\begin{eqnarray*}
	\|\uep(t)-\uep(s)\|_{2}\ \leq\
	\int_{s}^{t}\|\uep'(\tau)\|_{2}\,d\tau & \leq & |t-s|^{1/2}\left\{
	\int_{0}^{t}\|\uep'(\tau)\|_{2}^{2}\,d\tau\right\}^{1/2} \\
	 & \leq & |t-s|^{1/2}\left\{2\Fep(\uzep)\right\}^{1/2},
\end{eqnarray*}
so that the uniform bound on H\"{o}lder constants follows from 
assumption (\ref{hp:uzep}). 

Moreover, from (\ref{hp:uep-gf}) and statement~(\ref{stat:Linfty}) of 
Theorem~\ref{thm:uep} we have also that the functions 
$t\to\Fep(\uep(t))$ and $t\to\|\uep(t)\|_{\infty}$ are nonincreasing, 
hence
\begin{equation}
	\|\uep(t)\|_{\infty}+ \Fep(\uep(t))\leq
	\|\uzep\|_{\infty}+ \Fep(\uzep)
	\quad\quad
	\forall t\geq 0.
	\label{est:t-fix}
\end{equation}

Thanks to assumption~(\ref{hp:uzep}), the right-hand side is bounded
independently of $\ep$.  Therefore the compactness result in
statement~(2) of Theorem~\ref{thm:Gamma-conv} implies that the family
$\{\uep(t)\}$ is relatively compact in $L^{2}(\Omega)$ for every fixed
$t\geq 0$.

\paragraph{\emph{\textmd{Characterization of the limit}}}

We prove that, for every $T>0$, any limit point in the interval
$[0,T]$ of the family $\{\uep\}$ of approximating solutions is the
solution $u(t)$ of the total variation flow in the same interval.
Since the solution of the limit problem is unique, this is enough to
prove the convergence of the whole family.

Let $\ep_{n}\to 0^{+}$ be any sequence such that $u_{\ep_{n}}(t)$
uniformly converges to some $v(t)$ in $[0,T]$.  By
(\ref{hp:data-conv}) we have that $v(0)=u_{0}$.  So it is enough to
show that $v(t)$ is a maximal slope curve for the functional $TV(u)$
in $[0,T]$.

To this end, we apply Theorem~\ref{thm:main-cmp} to the sequence of
functionals $\{E^{**}_{\ep_{n}}(u)\}$.  Indeed they are convex
functionals, and they Gamma-converge to $TV(u)$ because of
Theorem~\ref{thm:Gamma-conv}.  Thus Remark~\ref{rmk:convex->scp} and
Proposition~\ref{prop:scp} prove that
assumption~(\ref{hp:slope-energy}) of Theorem~\ref{thm:main-cmp} is
satisfied.  Also (\ref{hp:sup-fn}) holds true because of
(\ref{est:t-fix}), and the fact that the functionals are nonnegative.

Since all the assumptions are satisfied, Theorem~\ref{thm:main-cmp}
implies that $v(t)$ is a maximal slope curve for the functional
$TV(u)$.

\paragraph{\emph{\textmd{Uniform convergence for all positive times}}}

It remains to prove that the convergence is global-in-time, as stated 
in (\ref{th:main}). This follows from two general facts.

The first one is that $u(t)$ tends, as $t\to +\infty$, to the constant
function $u_{\infty}$ equal to the average of $u_{0}$ in $\Omega$ 
(see~\cite{TV-property}).
The second fact is that $\uep(t)-u_{\infty}$ is the solution of
the approximating problem with
initial datum $u_{0\ep}-u_{\infty}$, hence 
statement~(\ref{stat:Linfty}) of Theorem~\ref{thm:uep} (with $p=2$) 
implies that
$t\to\|u_\ep(t)-u_{\infty}\|_{2}$ is a nonincreasing function.

Therefore for every 
$T>0$ we have that
\begin{eqnarray*}
	\sup_{t\geq T}\|\uep(t)-u(t)\|_{2} & \leq & 
	\sup_{t\geq T}\|\uep(t)-u_{\infty}\|_{2}+
	\sup_{t\geq T}\|u(t)-u_{\infty}\|_{2}\\
	 & = & \|\uep(T)-u_{\infty}\|_{2}+
	\sup_{t\geq T}\|u(t)-u_{\infty}\|_{2}  \\
	 & \leq & \|\uep(T)-u(T)\|_{2}+
	2\sup_{t\geq T}\|u(t)-u_{\infty}\|_{2},
\end{eqnarray*}
hence
\begin{eqnarray*}
	\sup_{t\geq 0}\|\uep(t)-u(t)\|_{2}& = & \max\Bigl\{
	\sup_{t\in[0,T]}\|\uep(t)-u(t)\|_{2},\; \sup_{t\geq
	T}\|\uep(t)-u(t)\|_{2}\Bigr\} \\
	 & \leq & \sup_{t\in[0,T]}\|\uep(t)-u(t)\|_{2}+
	 2\sup_{t\geq T}\|u(t)-u_{\infty}\|_{2}.
\end{eqnarray*}

Letting $\ep\to 0^{+}$, the first term tends to 0 because of the
convergence result in $[0,T]$.  Letting $T\to +\infty$, the second
term tends to 0 because $u(t)\to u_{\infty}$.  This proves
(\ref{th:main}) provided that initial data satisfy (\ref{hp:uzep}).
  
\paragraph{\emph{\textmd{Convergence for general data}}}

Let now $\uzep\to u_{0}$ be any family satisfying
(\ref{hp:data-conv}).  Let us choose a sequence $\{u_{0n}\}\subseteq
L^{\infty}(\Omega)\cap W^{1,\infty}(\Omega)$ with $u_{0n}\to u_{0}$.  For
every $n\in\n$, let $u_{\ep,n}(t)$ be the solution of the approximating
problem with $u_{\ep,n}(0)=u_{0n}$, and let $u_{n}(t)$ be the solution of
the limit problem with $u_{n}(0)=u_{0n}$.  For every $n\in\n$ we
already know that $u_{\ep,n}\to u_{n}$ in
$C^{0}\left([0,+\infty);L^{2}(\Omega)\right)$, because in this case 
the sequence of initial data does not depend on $\ep$ and satisfies 
(\ref{hp:uzep}).

Since both the approximating problems and the limit problem generate a 
contraction semigroup in $L^{2}(\Omega)$, we have that
\begin{eqnarray*}
	\|\uep(t)-u(t)\|_{2} & \leq & 
	\|\uep(t)-u_{\ep,n}(t)\|_{2}+\|u_{\ep,n}(t)-u_{n}(t)\|_{2}+
	\|u_{n}(t)-u(t)\|_{2}  \\
	 & \leq & \|\uzep-u_{0n}\|_{2}+\|u_{\ep,n}(t)-u_{n}(t)\|_{2}+
	\|u_{0n}-u_{0}\|_{2}.
\end{eqnarray*}

Taking the supremum over all $t\geq 0$, and letting $\ep\to 0^{+}$, 
we obtain that
$$\limsup_{\ep\to 0^{+}}\;\sup_{t\geq 0}\|\uep(t)-u(t)\|_{2}\leq
2\|u_{0n}-u_{0}\|_{2}.$$

Letting $n\to +\infty$, we finally obtain (\ref{th:main}) for general 
data.\qed

\subsubsection*{\centering Acknowledgments}

We would like to thank Patrick Guidotti for sending us a preliminary 
version of reference~\cite{patrick-DM}.

\label{NumeroPagine}

\end{document}